\documentclass[letterpaper, 10 pt, conference]{ieeeconf}

\IEEEoverridecommandlockouts
\usepackage{cite}

\usepackage{amsthm,amsmath,amssymb,amsfonts}
\usepackage{algorithmic}
\usepackage[ruled, noline, linesnumbered]{algorithm2e}
\usepackage{graphicx}
\usepackage{xcolor}

\usepackage{array}
\usepackage{multirow}

\makeatletter
\let\MYcaption\@makecaption
\makeatother

\usepackage[font=footnotesize]{subcaption}
\makeatletter
\let\@makecaption\MYcaption
\makeatother

\usepackage{textcomp}
\def\BibTeX{{\rm B\kern-.05em{\sc i\kern-.025em b}\kern-.08em
T\kern-.1667em\lower.7ex\hbox{E}\kern-.125emX}}

\usepackage{lipsum}

\newcommand\copyrighttext{
\footnotesize \copyright\ 2020 IEEE. Personal use of this material is permitted. Permission from IEEE must be obtained for all other uses, in any current or future media, including reprinting/republishing this material for advertising or promotional purposes, creating new collective works, for resale or redistribution to servers or lists, or reuse of any copyrighted component of this work in other works. DOI: 10.23919/ACC45564.2020.9147560
}
\newcommand\copyrightnotice{
\begin{tikzpicture}[remember picture,overlay]
\node[anchor=south,yshift=0.2in] at (current page.south) {\parbox{\dimexpr\textwidth-\fboxsep-\fboxrule\relax}{\copyrighttext}};
\end{tikzpicture}\vspace{-0.825\baselineskip}
}


\theoremstyle{plain}

\newtheorem{lemma}{Lemma}
\newtheorem{corollary}{Corollary}

\theoremstyle{definition}

\def\({\left(}
\def\){\right)}
\def\[{\left[}
\def\]{\right]}
\def\{{\left\lbrace}
\def\}{\right\rbrace}

\def\Embb{\mathbb{E}}  
\def\Rmbb{\mathbb{R}}  
\def\Ymbb{\mathbb{Y}}  
\def\Jcal{\mathcal{J}}  
\def\Mcal{\mathcal{M}}  
\def\Ncal{\mathcal{N}}  
\def\Rcal{\mathcal{R}}  
\def\Ucal{\mathcal{U}}  

\DeclareMathOperator*{\argmin}{argmin}

\newcommand{\E}[2][]{\Embb_{#1}\!\[#2\]}

\def\Comment#1{\hfill{\boldmath$\cdot$} #1\phantom{..}}
\def\Commentn#1{\hfill{\boldmath$\cdot$} #1\hspace{-0.5ex}}

\def\algcdc17{the algorithm in \cite{tseng2017local}}
\def\Algcdc17{The algorithm in \cite{tseng2017local}}
\def\sec#1{Section \ref{sec:#1}}
\def\fig#1{Fig. \ref{fig:#1}}
\def\subfig#1#2{Fig. \ref{fig:#1}\subref{subfig:#1-#2}}
\def\alg#1{Algorithm \ref{alg:#1}}
\def\lem#1{Lemma \ref{lem:#1}}
\def\cor#1{Corollary \ref{cor:#1}}
\def\eqn#1{\eqref{eq:#1}}

\newcommand{\centeritem}[3]{
\hspace*{#1}#2\hspace*{#3}

}

\newcommand{\centerlegend}[3]{\centeritem{#1}{\includegraphics{#2}}{#3}}

\title{\LARGE \bf A Generic Solver for Unconstrained Control Problems \\ with Integral Functional Objectives}

\author{Shih-Hao Tseng
\thanks{Shih-Hao Tseng is with the Division of Engineering and Applied Science, California Institute of Technology, Pasadena, CA 91125, USA. Email: {\tt\small shtseng@caltech.edu}}
\copyrightnotice
}

\usepackage{tikz}

\begin{document}

\maketitle
\thispagestyle{empty}
\pagestyle{empty}

\bstctlcite{IEEE_BSTcontrol}

\begin{abstract}
We present a generic solver for unconstrained control problems (UCPs) whose objectives take the form of an integral functional of the controllers. The solver generalizes and improves upon \algcdc17 for the Witsenhausen's counterexample, which provides the best-known results.
In essence, we show that minimizing the objective implies minimizing the marginal cost functions almost everywhere, and we perform the latter task pointwisely by the adaptive minimization technique, which speeds up the computation.
We implement single-threaded and parallelized versions of the proposed algorithm. Our implementation runs $30 \times$ faster than \algcdc17 on the Witsenhausen's counterexample, and we demonstrate the applicability of the solver and discuss the possible generalization to constrained problems and multi-dimensional controllers through three more examples.
\end{abstract}

\section{Introduction}\label{sec:introduction}
In this work, we focus on the unconstrained control problem (UCP) with the objective
\begin{align*}
\min\ \Jcal[ U ]
\end{align*}
where $U = \{u_0(y_0), u_1(y_1), \dots, u_{M-1}(y_{M-1})\}$ is the set of controllers $u_m : \Rmbb \to \Rmbb$. We assume that the objective functional $\Jcal[ U ]$ can be expressed as
\begin{align}
\Jcal[ U ] = \int L_m(u_m(y_m),y_m) \ d y_m + \Rcal_m[U_{-m}]
\label{eq:assumption-J}
\end{align}
for all $m = 0, \cdots, M-1$ (in other words, the functional can be expanded with respect to each $m$), where $L_m : \Rmbb \times \Rmbb \to \Rmbb$ and $\Rcal_m[U_{-m}]$ is the \emph{residual functional} depending only on the controllers other than $u_m$. UCP is ``unconstrained'' in the sense that we require the ability to embed all constraints, e.g., system dynamics, in \eqn{assumption-J}.
Although omitted, we remark that $L_m$ might also depend on $U_{-m}$ and the decomposition of $L_m$ and $\Rcal_m$ is not necessarily unique.

Finding the optimal controller of a control problem is a daunting task, even when the problem imposes no constraints.
The traditional approach to a control problem is to analyze its structure and conclude some useful properties that would help find the optimal controller. However, as demonstrated by the famous Witsenhausen's counterexample \cite{witsenhausen1968counterexample}, preferred properties, such as linearity, does not hold in general. As a result, deriving the optimal controller becomes craftsmanship relying on keen observations.

Fortunately, we often need a near-optimal controller rather than an exact optimal one in practice. There are two approaches, the analytical and the numerical, to design a near-optimal controller. The former approach examines some specific controller structure based on the problem characteristics. Again, the performance of the design highly depends on the sophisticated understanding of the problem. On the other hand, the numerical approach develops techniques to approximate the problem and obtain good approximations of the optimal controllers. Correspondingly, the main challenges lie on computation efficiency and approximation quality.

It used to be the computation-demanding nature of the numerical methods which impedes the adoption. But luckily, the advancing technologies in the past decades have hugely reshaped the research landscape: Cheaper computation resources and parallelization techniques facilitate the development of image classification \cite{krizhevsky2012imagenet}, machine learning \cite{goodfellow2014generative, abadi2016tensorflow}, and genomics processing \cite{turakhia2018darwin}. Increased computation power not only grants us higher efficiency but also enables finer sampling granularity and potentially better approximation quality. Therefore, we argue that numerical methods have great potential in the upcoming computation-rich era.

Although numerical methods could potentially compute faster with plenty of computation resources, its effectiveness plays a central role in getting closer to the optimum. Blindly adopting numerical methods can be ineffective.
Accordingly, it is critical to ask how to design numerical methods that are effective for general problems, in particular, for general UCPs.

Our approach to tackling general UCPs evolves from the algorithm proposed for the Witsenhausen's counterexample in 2017 \cite{tseng2017local}. \Algcdc17 outperforms all previous attempts on the well-known Witsenhausen's counterexample (e.g., \cite{li2009learning,karlsson2011iterative,mehmetoglu2014deterministic}), and its mechanism does not depend on the property of the given objective functional. Although the algorithm finds the controllers that result in the record-low cost, it is computationally demanding and requires a special math tool called calculus of variation.

\subsection{Contribution and Organization}

We examine UCP and provide a generic algorithm to find a near-optimal controller numerically. The proposed algorithm generalizes \algcdc17 with the following improvements:
First, it presents a new angle viewing the local Nash minimizing phase in \algcdc17 which does not involve calculus of variation.
In summary, our analysis reveals that UCP leads to a per-point marginal cost optimization problem, and the two methods used in \algcdc17, local Nash minimizing and local denoising, approach the problem using different candidate sets.
We also apply the adaptive minimization technique to speed up the convergence significantly. Furthermore, our proposed algorithm adopts a unified termination criterion applicable to arbitrary UCPs.

We then implement a generic solver based on the proposed algorithm in C++ and demonstrated that it converges $3 \times$ faster than \algcdc17 on Witsenhausen's counterexample. Since the proposed algorithm is parallelizable, we enhance the single-threaded C++ implementation using NVIDIA CUDA and observe another $10 \times$ computation speed up. We also demonstrate how the solver works on the zero-delay source-channel coding problem, inventory control problem, and $2$-dimensional Witsenhausen's counterexample. And we open source the tool for future research.

The paper is organized as follows. We first provide a brief overview of \algcdc17 in \sec{background}.
\sec{algorithm} introduces the ideas of marginal cost functions, local update, partial exhaustion, and adaptive minimization. Those ideas contribute to the design of our solver. We then implement the solver and examine its performance improvement in \sec{performance}. In \sec{examples}, we demonstrate that how we can use the solver to approach different problems. Finally, we conclude the paper in \sec{conclusion}.

\subsection{Notation}
By convention, we denote the state by $x$, control by $u$, observation by $y$, and disturbance by $w$. Let $\Embb_A$ be the expected value with respect to random variable $A$. We omit $A$ when the expectation is taken with respect to all random variables.
We denote by $A \sim \Ncal(\mu,\sigma^2)$ a Gaussian random variable $A$ with mean $\mu$ and variance $\sigma^2$ and by $A \sim \Ucal(a,b)$ a uniform random variable distributed over $[a,b]$.
Given a number $M$, we slightly abuse the notation to denote $m = 0, \dots, M-1$ by $m \in M$.
For a multivariate function $F(a,b)$, we introduce the shorthand notation $F'(a,b) = \frac{\partial F(a,b)}{\partial a}$ to denote the partial derivative with respect to the first variable.

Given a functional $\Jcal[U]$, we denote by $\frac{\delta \Jcal[U]}{\delta u}(y)$ the functional derivative of $\Jcal[U]$ with respect to the function $u(y)$, which is derived from the Taylor expansion below:
\begin{align*}
\Jcal[U + \epsilon \delta u] =&\ \Jcal[U] \\
+&\ \epsilon \int \frac{\delta \Jcal[U]}{\delta u}(y) \delta u(y) dy\\
+&\ \frac{\epsilon^2}{2} \int \frac{\partial}{\partial u} \frac{\delta \Jcal[U]}{\delta u}(y) \delta u^2(y) dy
+ O(\epsilon^3)
\end{align*}
where $U + \epsilon \delta u$ refers to the set $(U \backslash \{ u(y) \}) \cup \{ u(y) + \epsilon \delta u(y) \}$, $\delta u(y)$ is a bounded function for variation, and $O(\epsilon^3)$ represents the residual terms of order $\epsilon^3$ or higher.

\section{Existing Algorithm and its Limitation}\label{sec:background}

We first explain \algcdc17, which attains the state-of-the-art best results for Witsenhausen's counterexample, in \sec{background-basic-structure}. We then discuss the limitations of the algorithm in \sec{background-limitations}.

\subsection{Basic Structure}\label{sec:background-basic-structure}
In \cite{tseng2017local}, Witsenhausen's counterexample is deemed an optimization problem minimizing a given functional. To obtain the best controllers, \algcdc17 introduces two main components: local Nash minimizers and local denoising.

\subsubsection{Local Nash Minimizers}
\cite{tseng2017local} shows that an optimal controller must be a local Nash minimizer. Using calculus of variation, the first order condition (FOC) and the second order condition (SOC)
\begin{align*}
\frac{\delta \Jcal[U]}{\delta u_m}(y_m) = 0,\quad \frac{\partial}{\partial u_m}\frac{\delta \Jcal[U]}{\delta u_m}(y_m) \geq 0
\end{align*}
are then derived for local Nash minimizers. Combining FOC and SOC, \algcdc17 repeats revised Newton's method to seek for a local Nash minimizer.

\subsubsection{Local Denoising}
The most important observation given by \cite{tseng2017local} is that finding local Nash minimizers numerically would land in a ``noisy'' controller. As a result, \cite{tseng2017local} introduces the idea of ``denoising,'' i.e., for each $u_m$, we denoise for all $y_m$ by
\begin{align*}
u_m(y_m) \leftarrow&\ \argmin\limits_{u_m(y) : y \in B_r(y_m)} C_m(u_m(y),y_m),
\end{align*}
where $B_r(a)$ is a ball centering at $a$ with radius $r$.

\subsection{Limitations}\label{sec:background-limitations}
Although \algcdc17 is demonstrated effective for Witsenhausen's counterexample and it is applicable to other similar problems such as inventory control, there are still few issues left by \cite{tseng2017local}. First, albeit a universal approach to finding local Nash minimizers, calculus of variation is not a simple idea/operation for the people who know little about functional analysis. For local denoising, \cite{tseng2017local} obtains functions $C_m$ by observation. It would be more rigorous to have a standard procedure to derive $C_m$. Meanwhile, despite its great performance in terms of the final cost it achieves, the algorithm runs slow in practice. And it is not clear how the termination criterion used in \cite{tseng2017local} can be easily generalized for arbitrary problems.

\section{Generic Algorithm Design}\label{sec:algorithm}
In this section, we study UCP and illustrate how to overcome the limitations stated in \sec{background-limitations} so that we can improve the ideas in \cite{tseng2017local} to solve UCPs.

\subsection{Marginal Cost Functions}
We start the analysis with \lem{optimal-criterion}, which is a necessary condition for optimal $U$, followed by the derivation of the \emph{marginal cost function} $C_m$.

\begin{lemma}\label{lem:optimal-criterion}
If $U$ minimizes $\Jcal$, we have
\begin{align*}
u_m(y_m) = \argmin\limits_{u \in \Rmbb} L_m(u,y_m)
\end{align*}
almost everywhere for all $m \in M$.
\end{lemma}

The lemma can be derived from \eqn{assumption-J} and the derivation is straightforward. Suppose \lem{optimal-criterion} is not true, there must exist some $\Delta > 0$ such that
\begin{align*}
L_m(u_m(y_m),y_m) \geq \Delta + \min\limits_{u \in \Rmbb} L_m(u,y_m)
\end{align*}
over some set $Y_{\Delta}$ with non-zero measure $\Mcal(Y_{\Delta})$. As such, we can set $u_m(y_m)$ as in \lem{optimal-criterion} for all $y_m \in Y_{\Delta}$, which results in the reduction of the functional value by at least $\Mcal(Y_{\Delta}) \Delta > 0$. However, it leads to a contradiction as $U$ minimizes $\Jcal$.

\lem{optimal-criterion} relates the optimal $u_m(y_m)$ with $L_m(u,y_m)$. However, a UCP is usually specified by $\Jcal$ and its decomposition to $L_m(u,y_m)$ is in general not unique. We would prefer to base our solver on \lem{optimal-criterion} with respect to a more deterministic expression than $L_m$. Therefore, we derive the \emph{marginal cost function} $C_m$ as follows.

From \eqref{eq:assumption-J}, we know
\begin{align*}
L_m(u_m(y_m),y_m) = \frac{\partial \Jcal[U]}{\partial y_m} - \frac{\partial \Rcal_m[U_{-m}]}{\partial y_m}
\end{align*}
which implies that
\begin{align*}
L_m(u,y_m) = \left. \( \frac{\partial \Jcal[U]}{\partial y_m} - \frac{\partial \Rcal_m[U_{-m}]}{\partial y_m} \) \right|_{u_m(y_m) = u}.
\end{align*}
Substitute it into \lem{optimal-criterion}, we have
\begin{align*}
u_m(y_m)
=&\ \argmin\limits_{u \in \Rmbb} \left. \( \frac{\partial \Jcal[U]}{\partial y_m} - \frac{\partial \Rcal_m[U_{-m}]}{\partial y_m} \) \right|_{u_m(y_m) = u} \\
=&\ \argmin\limits_{u \in \Rmbb} \left. \frac{\partial \Jcal[U]}{\partial y_m} \right|_{u_m(y_m) = u}
\end{align*}
because $\Rcal_m[U_{-m}]$ does not depend on $u$.
Thus, we define the \emph{marginal cost function} $C_m$ for each $m \in M$ by
\begin{align*}
C_m(u,y_m) = \left. \frac{\partial \Jcal[U]}{\partial y_m} \right|_{u_m(y_m) = u}
\end{align*}
and we can rephrase \lem{optimal-criterion} as
\begin{corollary}\label{cor:optimal-criterion}
If $U$ minimizes $\Jcal$, we have
\begin{align*}
u_m(y_m) = \argmin\limits_{u \in \Rmbb} C_m(u,y_m)
\end{align*}
almost everywhere for all $m \in M$.
\end{corollary}

Now, once we have the objective $\Jcal$ specified, $C_m$ can be derived using simple partial derivation. Actually, the definition of $C_m$ also matches with the observation in \cite{tseng2017local}.

\subsection{Local Update and Partial Exhaustion}\label{sec:algorithm-local-update-partial-exhaustion}
\cor{optimal-criterion} is not only a necessary constraint but also a way to improve a non-optimal controller. Essentially, we can always improve a non-optimal solution by
\begin{align*}
u_m(y_m) \leftarrow \argmin\limits_{u \in \Rmbb} C_m(u,y_m).
\end{align*}
To do so, we need to find the minimizer of $C_m(u,y_m)$.

To find the minimizer, the most effective method is to exhaust all possible $u \in \Rmbb$ and choose the best one. A full exhaustion is rarely practical as the computational cost would be intolerable. Instead, we would prefer computationally economic methods, such as \emph{local update} and \emph{partial exhaustion}.

Local update methods include gradient descent, Newton's method, and their variations. Those methods rely on local information at $y_m$, such as derivatives, to update $u_m$. Partial exhaustion takes a different approach. It avoids the heavy computational load in a full exhaustion by searching within only a subset of candidate values.

Not surprisingly, these two methods correspond to the two main components in \algcdc17.
Since
\begin{align*}
\frac{\delta \Jcal[U]}{\delta u_m}(y_m) =&\ C_m'(u_m(y_m),y_m), \\
\frac{\partial}{\partial u_m}\frac{\delta \Jcal[U]}{\delta u_m}(y_m) =&\ C_m''(u_m(y_m),y_m),
\end{align*}
finding local Nash minimizers using the revised Newton's method in \cite{tseng2017local} is equivalent to minimizing the function $C_m(u,y_m)$ at each $y_m$ by a revised Newton's method.

\begin{figure}[!t]
\centering
\subcaptionbox{Local update\label{subfig:candidate-set-local-update}}[0.48\columnwidth][c]{
\includegraphics[scale=0.5]{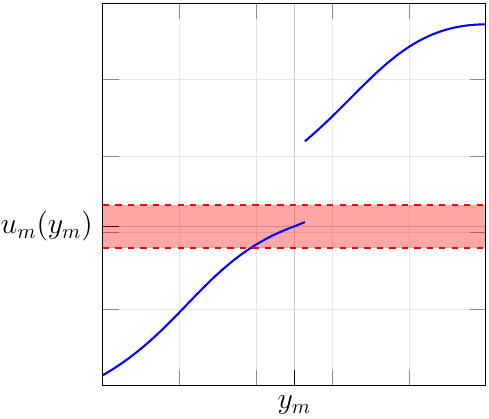}
{\footnotesize $\{ u : u \in B_r(u_m(y_m)) \}$}
}\hfill
\subcaptionbox{Partial exhaustion\label{subfig:candidate-set-partial-exhaustion}}[0.48\columnwidth][c]{
\includegraphics[scale=0.5]{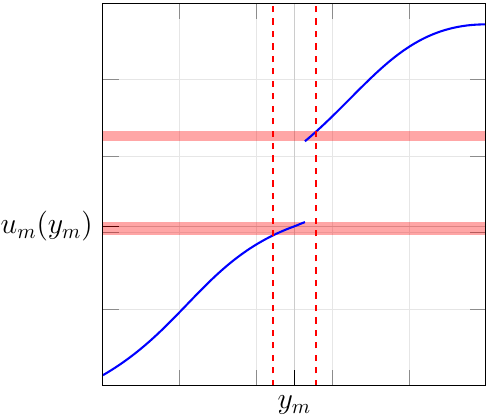}
{\footnotesize $\{ u_m(y) : y \in B_r(y_m) \}$}
}
\caption{The candidate sets used in different techniques are marked by shaded areas. We depict the function $u_m$ by the thick blue line and mark the range of $B_r$ by red dashed lines. The candidate set can be disjoint in partial exhaustion when $u_m(y_m)$ is discontinuous.}
\label{fig:candidate-set}
\end{figure}

On the other hand, local denoising is a partial exhaustion method that searches within the candidate set
\begin{align}
\{ u_m(y) : y \in B_r(y_m) \}
\label{eq:partial-exhaustion-candidate-set}
\end{align}
to minimize $C_m(u,y_m)$. This is an important feature of \algcdc17: It searches over the domain instead of the range of the function $u_m$. Usually, a partial exhaustion method searches the argument value locally, i.e., it considers the candidate set
\begin{align*}
\{ u : u \in B_r(u_m(y_m)) \}.
\end{align*}
Such a partial exhaustion method plays a similar role as a local update method since they both try to update
\begin{align*}
u_m(y_m) \leftarrow&\ \argmin\limits_{u : u \in B_r(u_m(y_m))} C_m(u,y_m)
\end{align*}
for some small $r$. As a result, we do not benefit much from combining the methods together.
However, \eqref{eq:partial-exhaustion-candidate-set} can be a disconnected set when $u_m$ is noncontinuous as in \fig{candidate-set}, which allows searching and ``leaping'' to another region.

\subsection{Adaptive Minimization}\label{sec:algorithm-adaptive-minimization}
Both local update and partial exhaustion methods are adopted alternatively in \algcdc17: In each repetition round, it denoises the controllers after repeating the revised Newton's method several times. This hybrid strategy converges to the best-known results, but it progresses quite slowly. The reason is that the two methods improve $\Jcal$ in different ways and one may be more effective than the other at different searching phase.

\begin{figure}[!t]
\centering
\centerlegend{0.445in}{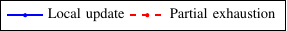}{0.07in}
\subcaptionbox{\Algcdc17 runs local update $19$ times and performs $1$ partial exhaustion per round. Partial exhaustion is more effective than local update.\label{subfig:improvements-cdc2017}}{
\includegraphics[scale=0.5]{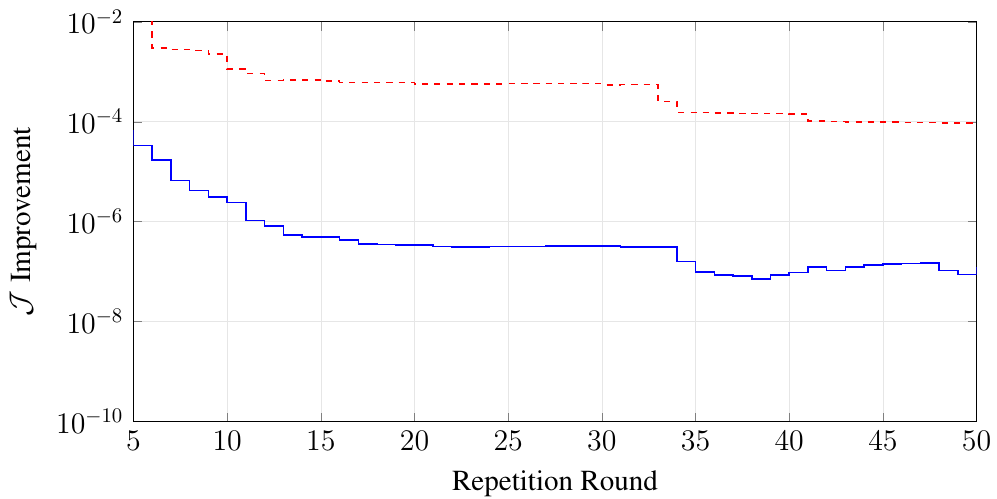}
}
\vspace*{0.5\baselineskip}

\centerlegend{0.445in}{method_legends}{0.07in}
\subcaptionbox{Adaptive minimization allocates $20$ iterations in each round to the methods according to the improvement in the last round. As such, both methods improve $\Jcal$ comparably. \label{subfig:improvements-parallel}}{
\includegraphics[scale=0.5]{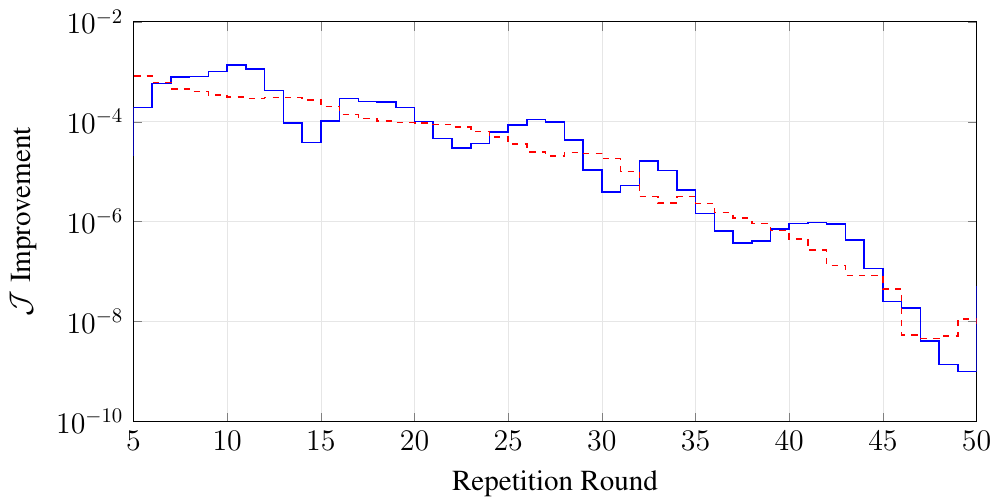}
}
\vspace*{0.5\baselineskip}

\centerlegend{0.355in}{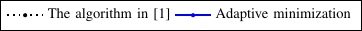}{0.055in}
\vspace*{0.1in}

\subcaptionbox{Comparison of the overall improvement. Adaptive minimization outperforms \algcdc17 significantly.\label{subfig:improvements-comparison}}{
\includegraphics[scale=0.5]{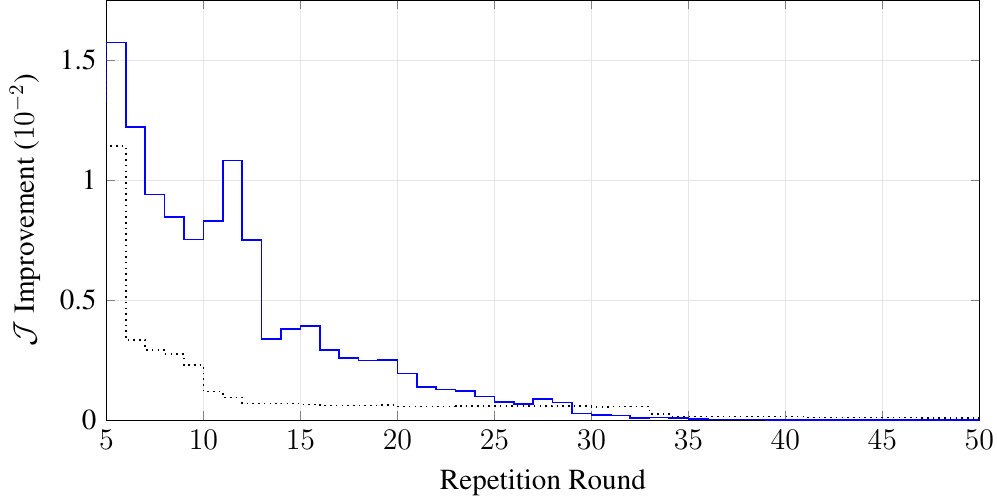}
}

\caption{The improvement of the objective functional $\Jcal$ in the Witsenhausen's counterexample for each method.
Each round has $20$ iterations. }
\label{fig:improvements}
\end{figure}

In \subfig{improvements}{cdc2017}, we consider the Witsenhausen's counterexample as in \cite{tseng2017local} and plot the average $\Jcal$ value improvement of the local update method (revised Newton's method) and the partial exhaustion method (denoising) for a series of repetition rounds, each round consisting of $19$ local update iterations and $1$ partial exhaustion iteration. \subfig{improvements}{cdc2017} shows that partial exhaustion performs much better than local update, but local update runs much more times than partial exhaustion. As a result, the overall improvement at each round is low as in \subfig{improvements}{comparison}.

To improve the efficiency, we introduce \emph{adaptive minimization}. The basic idea is that we will adapt the number of iterations according to the performance of the method. The more effective method gets more iterations to run. Fixing the total iterations to be $20$ per round, we perform adaptive minimization and \subfig{improvements}{parallel} shows that the average improvements of the two methods are comparable. Furthermore, the overall improvement is boosted by adaptive minimization as shown in \subfig{improvements}{comparison}.

\subsection{Generic Algorithm}
We combine the methods described in \sec{algorithm-local-update-partial-exhaustion} and \sec{algorithm-adaptive-minimization} to construct \alg{generic-solver} for UCP. In summary, \alg{generic-solver} approximates the controllers $U$ by step functions over a given sample range. Then, it adaptively minimizes $\Jcal[U]$ using local update and partial exhaustion methods. Adaptive minimization is repeated until some precision factor is met. The whole procedure can be partitioned into five main parts, and we give more detailed descriptions below.

\begin{algorithm}[!t]
\begin{algorithmic}[1]
\REQUIRE Number of iterations per round $N$ and\\ \hspace*{2.5ex} precision $p \geq 0$.
\STATE{Initialize controllers $U$.}\label{alg:generic-solver:initialize-begin}
\STATE{$\Jcal_c \leftarrow \Jcal[U]$.
\Comment{Initial objective value}}
\STATE{$I_L \leftarrow p$,
$I_P \leftarrow p$.
\Comment{Initial improvements}}
\STATE{$N_L \leftarrow \left\lfloor \frac{N}{2} \right\rfloor$,
$N_P \leftarrow N - N_L$.
\Comment{Initial iterations}}\label{alg:generic-solver:initialize-end}
\makeatletter\ALC@it\algorithmicwhile\ {$I_L + I_P > p$}\ \algorithmicdo
\Commentn{Main loop}
\begin{ALC@whl}\label{alg:generic-solver:main-loop}
\FOR{$N_L$ iterations}\label{alg:generic-solver:local-update-begin}
\FOR{$m = 0$ \TO $M-1$}
\STATE{LocalUpdate($u_m$)}
\ENDFOR
\ENDFOR
\STATE{$\Jcal_c \leftarrow \Jcal[U]$,
$I_L \leftarrow |I_L - \Jcal_c|$.\\
\Comment{Get local update improvement}}\label{alg:generic-solver:local-update-end}
\FOR{$N_P$ iterations}\label{alg:generic-solver:partial-exhaustion-begin}
\FOR{$m = 0$ \TO $M-1$}
\STATE{PartialExhaustion($u_m$)}
\ENDFOR
\ENDFOR
\STATE{$\Jcal_c \leftarrow \Jcal[U]$,
$I_P \leftarrow |I_P - \Jcal_c|$.\\
\Comment{Get partial exhaustion improvement}}\label{alg:generic-solver:partial-exhaustion-end}
\STATE{$N_L = \min\{\max\{ \left\lfloor \frac{I_L N}{I_L + I_P} \right\rfloor, 1\}, N-1\}$,
$N_P \leftarrow N - N_L$.
\Comment{Adaptive minimization}}\label{alg:generic-solver:adaptive-minimization}
\ENDWHILE
\end{algorithmic}
\caption{Generic Solver}
\label{alg:generic-solver}
\end{algorithm}

\vspace*{0.3\baselineskip}
\subsubsection{Initialization (line \ref{alg:generic-solver:initialize-begin} -- \ref{alg:generic-solver:initialize-end})}
Given some sampling range $[a,b] \subseteq \Rmbb$ and number of samples $d$, we create a linearly spaced vector $\Ymbb \in \Rmbb^d$ spanning over $[a,b]$ as the domain. Each controller $u_m$ is then approximated by a function mapping $\Ymbb$ to $\Rmbb$. As suggested in \cite{tseng2017local}, we initialize $u_m$ by an identity map:
\begin{align*}
u_m(y_m) \leftarrow y_m \quad \text{for all} \quad y_m \in \Ymbb.
\end{align*}

Also, we initialize the current objective value $\Jcal_c$, improvement $I_L, I_P$, and iterations $N_L, N_P$ for adaptive minimization. We use subscript $L$ to denote the variable for local update and subscript $P$ for partial exhaustion.

\vspace*{0.3\baselineskip}
\subsubsection{Main Loop (line \ref{alg:generic-solver:main-loop})}
We repeat adaptive minimization in the main loop, and we adopt a unified termination criterion: The loop terminates when the overall improvement in the current round is smaller than the precision $p$. This criterion improves upon the one adopted in \cite{tseng2017local} as it is valid regardless of the objective $\Jcal$.

\vspace*{0.3\baselineskip}
\subsubsection{Local Update (line \ref{alg:generic-solver:local-update-begin} -- \ref{alg:generic-solver:local-update-end})}
In this phase, we aim to solve
\begin{align*}
u_m(y_m) \leftarrow&\ \argmin\limits_{u : u \in B_r(u_m(y_m))} C_m(u,y_m)
\end{align*}
for all $y_m \in \Ymbb$ using some local update methods.

In our implementation, we adopt a modified Newton's method slightly different from the one in \cite{tseng2017local}: We only apply Newton's method
\begin{align*}
u_m(y_m) \leftarrow u_m(y_m) - \frac{C_m'(u_m(y_m),y_m)}{C_m''(u_m(y_m),y_m)}
\end{align*}
when
\begin{align*}
C_m''(u_m(y_m),y_m) > 0
\end{align*}
which implies a local minimum. Otherwise, we apply simple gradient method
\begin{align*}
u_m(y_m) \leftarrow u_m(y_m) - \tau C_m'(u_m(y_m),y_m)
\end{align*}
where $\tau$ is the given step size.

\vspace*{0.3\baselineskip}
\subsubsection{Partial Exhaustion (line \ref{alg:generic-solver:partial-exhaustion-begin} -- \ref{alg:generic-solver:partial-exhaustion-end})}
We perform partial exhaustion based on the sampled candidate set \eqref{eq:partial-exhaustion-candidate-set}:
\begin{align*}
u_m(y_m) \leftarrow&\ \argmin\limits_{u_m(y) : y \in B_r(y_m) \cap \Ymbb} C_m(u_m(y),y_m)
\end{align*}
for all $y_m \in \Ymbb$. Essentially, it is equivalent to the local denoising procedure with denoising radius $r$. The radius $r$ should be chosen such that it allows the algorithm to examine at least $1$ sample point each direction along $\Rmbb$. Meanwhile, $r$ should not be too large as the larger $r$ leads to higher computation load and potentially traps the local update phase around suboptimal local minima. In our implementation, we choose $r$ such that we check $2$ points per direction.

\vspace*{0.3\baselineskip}
\subsubsection{Adaptive Minimization (line \ref{alg:generic-solver:adaptive-minimization})}
We perform a per-round adaptive minimization by updating the iterations according to the improvements in the last round. Starting with fair sharing as in line \ref{alg:generic-solver:local-update-end}, we calculate the improvements gained from each of the methods by line \ref{alg:generic-solver:local-update-end} and \ref{alg:generic-solver:partial-exhaustion-end}.
In line \ref{alg:generic-solver:adaptive-minimization}, we allocate the $N$ iterations in the next round proportional to the contribution of each method while ensuring each method will be run at least once.

\section{Performance}\label{sec:performance}
To demonstrate the performance, we implement \alg{generic-solver} in two different versions: the single-threaded solver in C++ (denoted by \emph{Single}) and the parallelized one in NVIDIA CUDA (denoted by \emph{Parallel}). The solvers are open sourced at \cite{UCP_solver}. The solvers are run under the number of iterations per round $N = 20$ and the precision $p = 10^{-10}$.

We compare our solver with \algcdc17
on the Witsenhausen's counterexample, for which the formulation is deferred to \sec{examples-witsenhausen}, under the parameters $k = 0.2$ and $\sigma = 5$.

For fair comparison, we
impose the same termination criterion and total iterations per round for \algcdc17: It performs local Nash minimizing $19$ times followed by $1$ local denoising.
All experiments are run on a desktop with Intel Xeon W-2145 (CPU), NVIDIA Quadro P1000 (GPU), and Samsung 970 EVO NVMe M.2 (SSD).

\subsection{Parallelization}
To reduce the computation time, we examine \alg{generic-solver} and identify the parts that can be parallelized.

In both the local update and partial exhaustion phases, we revise $u_m$ one at a time, which may change $C_{m'}$ for some $m' \in M, m' \neq m$. As a result, such a dependency prevents the computation being done in parallel. On the other hand, updating $u_m(y_m)$ does not affect $C_m$, and hence the calculation for each $y_m \in \Ymbb$ can be parallelized to reduce the computation time.

\subsection{Computation Time}

\begin{figure}[!t]
\centering
\centerlegend{0.485in}{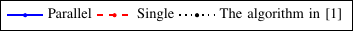}{0.075in}
\includegraphics[scale=0.5]{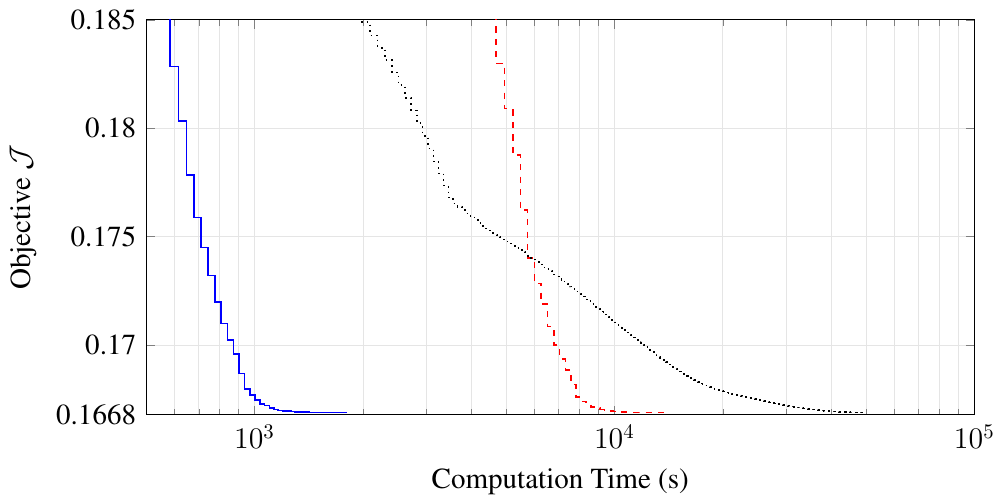}
\caption{Comparison of the convergence time among different methods. \algcdc17 converges faster in the beginning but it converges much slower when approaching to the optimum. On the contrary, both single-threaded and parallelized versions of \alg{generic-solver} converge much faster towards the optimum.}
\label{fig:time-performance}
\end{figure}

\fig{time-performance} shows how fast the solvers and \algcdc17 converges along the computation time under $14000$ sample points. \algcdc17 converges quickly at the beginning since it computes $u_1$ directly through a closed form expression. On the contrary, \alg{generic-solver} does not utilize any closed form expression. It simply searches for $u_0$ and $u_1$ using local update and partial exhaustion. Although \algcdc17 is more effective in computing $u_1$, adaptive minimization allows \alg{generic-solver} to select and use the more effective method. As a result, the single-threaded version converges slower at the beginning, but it improves $\Jcal$ much more effectively than \algcdc17 when getting closer to the limit. After parallelization, the convergence speed is further improved by $10 \times$, shifting the curve to the left in \fig{time-performance}.

\subsection{Scalability}
\begin{figure}[!t]
\centering
\centerlegend{0.355in}{legends}{0.055in}
\includegraphics[scale=0.5]{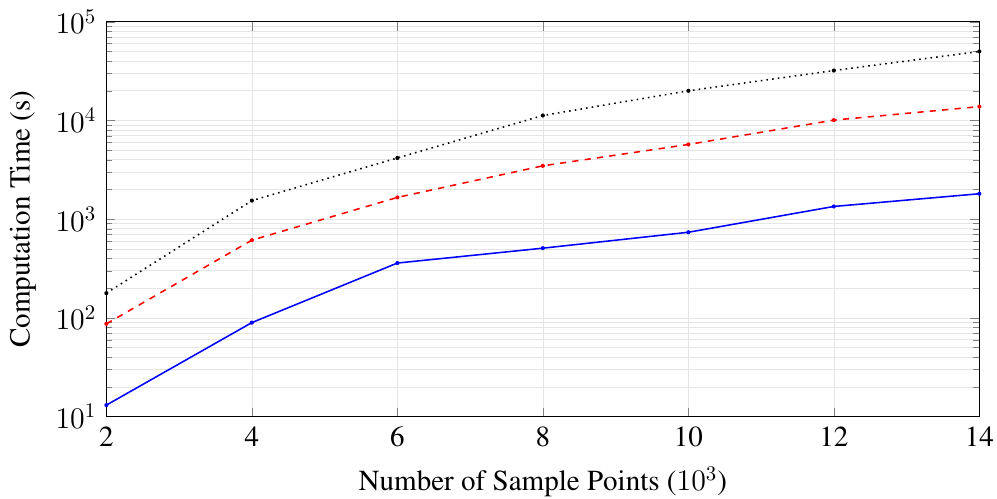}
\caption{\alg{generic-solver} scales much better than \algcdc17. The parallelized version generates the controllers $30 \times$ faster than \algcdc17.}
\label{fig:scalability}
\end{figure}

Another property we examine is scalability. Scalability is important for a numerical method as it implies how precise the approximation can be and how large the problem the solver can handle. To demonstrate the scalability, we vary the number of sample points for each controller and measure the computation time needed before reaching the required precision level. In \fig{scalability}, the single-threaded solver runs about $3 \times$ faster than \algcdc17.
After parallelization, the computation time is further reduced by $10 \times$, and hence it enjoys $30 \times$ performance improvement over \algcdc17.

\section{Examples}\label{sec:examples}
We apply our solver to three different examples: Witsenhausen's counterexample (\sec{examples-witsenhausen}), zero-delay source-channel coding (\sec{examples-coding}), and inventory control (\sec{examples-inventory}). Since inventory control problem imposes additional constraints on the controllers, we discuss how our solver might be generated for the problem with constraints.

\def\figsep{\hspace*{1.5ex}}

\subsection{Witsenhausen's Counterexample}\label{sec:examples-witsenhausen}

\begin{figure}[!t]
\centering
\includegraphics[scale=0.5]{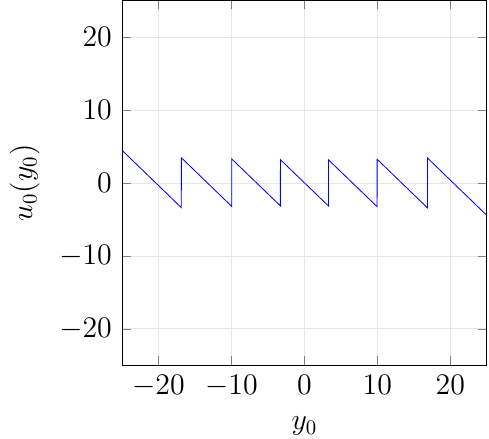}\figsep
\includegraphics[scale=0.5]{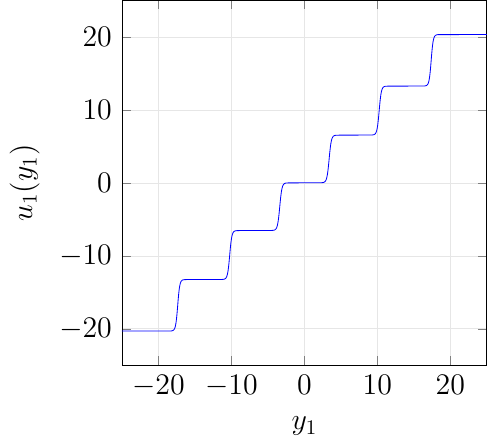}
\vspace*{-0.5\baselineskip}
\caption{Controllers for Witsenhausen's counterexample ($k = 0.2, \sigma = 5$),
$\Jcal[U] = 0.166897$.
}
\label{fig:example-witsenhausen}
\end{figure}

Witsenhausen's counterexample \cite{witsenhausen1968counterexample} aims to minimize the objective
\begin{align*}
\min\ \Jcal[U] = \E{ k^2 u_0(y_0)^2 + x_2^2 }
\end{align*}
subject to the system dynamic
\begin{alignat*}{2}
x_1 =&\ x_0 + u_0(y_0), &\ \quad \ y_0 =&\ x_0, \nonumber \\
x_2 =&\ x_1 - u_1(y_1), &\ \quad \ y_1 =&\ x_1 + w,
\end{alignat*}
where $x_0 \sim \Ncal(0, \sigma^2)$ and $w \sim \Ncal(0,1)$.

We use our solver to find controllers in \fig{example-witsenhausen}. The result is very close to the best-known controller in \cite{tseng2017local}.

\subsection{Zero-Delay Source-Channel Coding}\label{sec:examples-coding}

\begin{figure}[!t]
\centering
\includegraphics[scale=0.5]{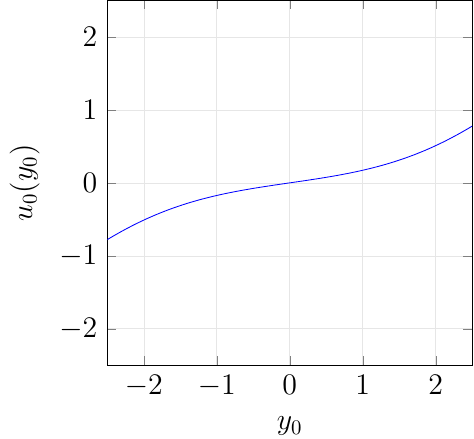}\figsep
\includegraphics[scale=0.5]{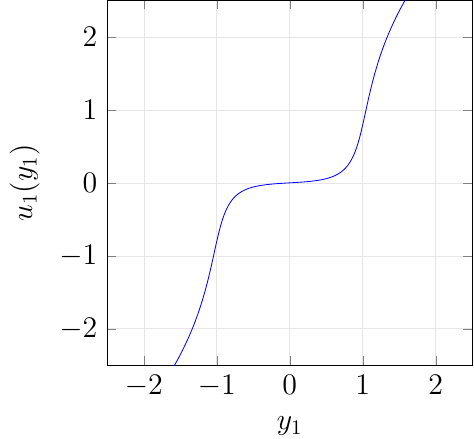}
\vspace*{-0.5\baselineskip}
\caption{Controllers for zero-delay source-channel coding problem ($\lambda = 2$),
$\Jcal[U] = 0.890756$.}
\label{fig:example-coding}
\end{figure}

The zero-delay source-channel coding problem \cite{akyol2014zero} has the objective
\begin{align*}
\min\ \Jcal[U] = \E{ \lambda u_0(x_0)^2 + (u_1(x_1) - x_0)^2 }
\end{align*}
and the system dynamic
\begin{align*}
y_0 = x_0, \quad x_1 = u_0(y_0) + w, \quad y_1 = x_1,
\end{align*}
where $x_0 \sim \Ncal(0, 1)$ and $w \sim \Ucal(-1,1)$.

\fig{example-coding} shows the non-linear controllers found by our solver.

\subsection{Inventory Control and Constrained Controller}\label{sec:examples-inventory}
\begin{figure}[!t]
\centering
\includegraphics[scale=0.5]{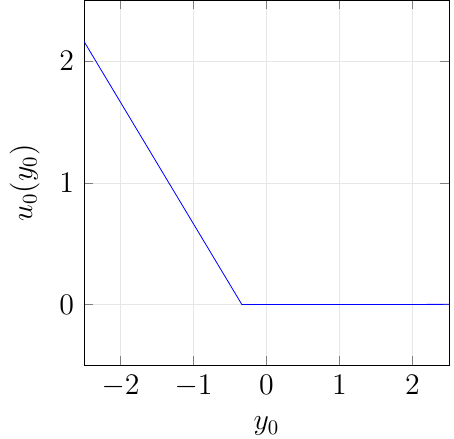}\figsep
\includegraphics[scale=0.5]{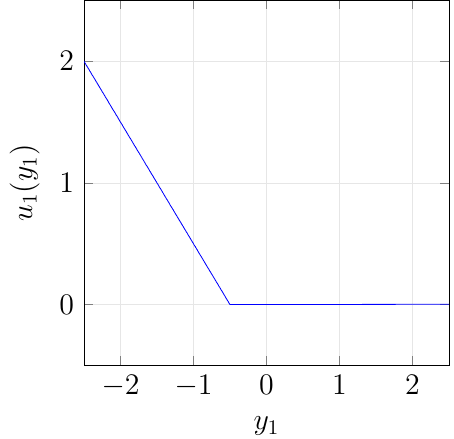}
\vspace*{-0.5\baselineskip}
\caption{Controllers for the inventory control problem using the same setting as in \cite{tseng2017local}.}
\label{fig:example-inventory}
\end{figure}

The inventory control problem has the objective
\begin{align*}
\min\ \Jcal[U] = \E{ \sum\limits_{m=0}^M \xi u_m(x_m) +  \gamma(x_{m+1}) }
\end{align*}
where $\gamma$ is some cost function attaining the minimum at $0$.

The system dynamic is
\begin{align*}
y_m = x_m, \quad x_{m+1} = x_m + u_m(x_m) - w_m,
\end{align*}
where $x_0 \sim \Ucal(-1,1)$ and $w_m \sim \Ucal(-1,1)$ for all $m \in M$.

The inventory control problem imposes an additional constraint: $u_m \geq 0$. We can enforce the constraint by performing $u_m \leftarrow \max\{ 0 , u_m\}$ after each local update and partial exhaustion. With the same settings as in \cite{tseng2017local}, our solver finds the same controllers as in \fig{example-inventory}. It is reported in \cite{tseng2017local} that this enforcement successfully finds the optimal controllers. Therefore, we conjecture that it is possible to generalize this generic UCP solver for the problems with constraints by projecting the result back to the feasible region after each local update and partial exhaustion.

\subsection{Multi-Dimensional Controller: $2$-Dimensional Witsenhausen's Counterexample}
We can also apply the solver to multi-dimensional controllers. For example, $2$-dimensional Witsenhausen's counterexample \cite{grover2013approximately,subramanian2018some} aims to minimize the objective
\begin{align*}
\min\ \Jcal[U] = \frac{1}{2}\E{ k^2 \lVert u_0(y_0) \rVert^2 + \lVert x_2 \rVert^2 }
\end{align*}
under the same state dynamic as the Witsenhausen's counterexample \cite{witsenhausen1968counterexample}. The difference is that the state $x$, controller $u$, and output $y$ are all $2$-dimensional vectors.

\begin{figure}[!t]
\centering
\centeritem{0.301in}{\includegraphics[scale=0.5]{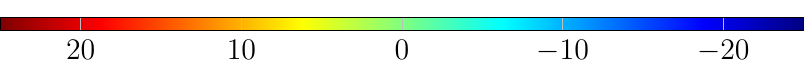}}{0.044in}

\includegraphics[scale=0.5]{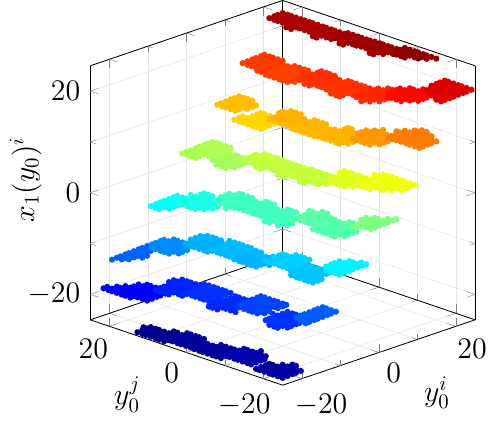}\figsep
\includegraphics[scale=0.5]{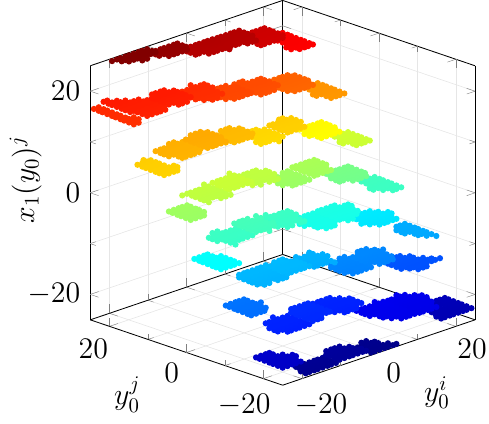}
\vspace*{0.25\baselineskip}

\includegraphics[scale=0.5]{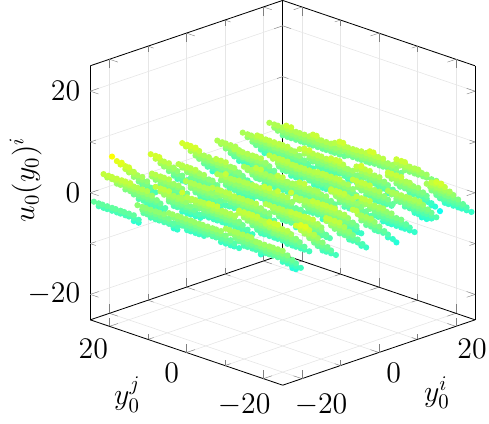}\figsep
\includegraphics[scale=0.5]{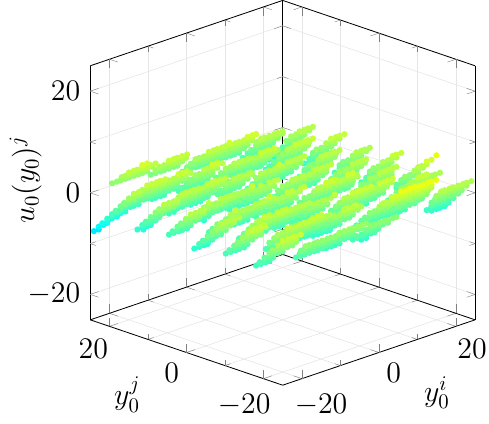}
\vspace*{0.25\baselineskip}

\includegraphics[scale=0.5]{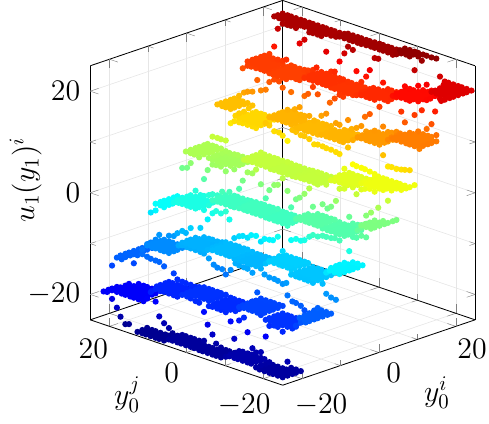}\figsep
\includegraphics[scale=0.5]{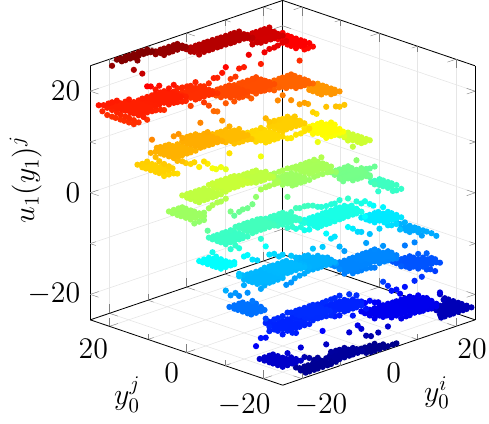}

\vspace*{-0.5\baselineskip}
\caption{State $x_1$ and controllers $u_0$ and $u_1$ for $2$-dimensional Witsenhausen's counterexample ($k = 0.2, \sigma = 5$). We denote the two dimensions by superscripts $i$ and $j$. With only $25600$ sample points over $y_0$ and $y_1$, we can obtain the result $\Jcal[U] = 0.166719$, which is close to the best known one $0.1527$ by machine learning and sophisticated heuristics in \cite{subramanian2018some}.
}
\label{fig:example-witsenhausen-2d}
\end{figure}

We apply our solver with only $160 \times 160 = 25600$ sample points over $y_0$ and $y_1$ and obtain the results in \fig{example-witsenhausen-2d}. The resulting objective value is very close to the best-known one in \cite{subramanian2018some}.

\section{Conclusion}\label{sec:conclusion}

The paper examines the unconstrained control problems and proposes an effective generic solver which can serve as the benchmark for the future UCP research. On the other hand, a lot of control problems do impose constraints on either the states of the controllers.
As a result, it would be of interest to generalize the solver to constrained control problems.

\bibliographystyle{IEEEtran}
\bibliography{Test}

\end{document}